\documentclass[11pt]{article}
\usepackage{epsfig}

\def\be{\begin{equation}}
\def\ee{\end{equation}}
\def\bea{\begin{eqnarray}}
\def\eea{\end{eqnarray}}
\def\bes{\begin{eqnarray*}}
\def\ees{\end{eqnarray*}}

\def\nn{\nonumber}
\def\lb{\label}
\def\bs{\setminus}
\def\pt{\partial}

\def\R{{\bf R}}
\def\C{{\bf C}}
\def\Z{{\bf Z}}
\def\K{{\bf K}}
\def\N{{\bf N}}
\def\U{{\bf U}}

\def\Q{{\bf Q}}

\def\aa{{\alpha}}
\def\bb{{\beta}}
\def\ga{{\gamma}}

\def\th{{\theta}}
\def\om{{\omega}}
\def\Om{{\Omega}}

\def\lm{{\lambda}}
\def\Lm{{\Lambda}}

\def\sg{{\sigma}}
\def\dm{{\diamond}}

\def\vf{{\varphi}}

\def\K{{\cal K}}
\def\P{{\cal P}}

\def\rank{{\rm rank}}
\def\Sp{{\rm Sp}}

\def\mod{{\rm mod}}
\def\ol{\overline}

\def\hb{\vrule height0.18cm width0.14cm $\,$}

\title{On the multiplicity of closed geodesics on Finsler spheres with irrationally elliptic closed geodesics}

\author{Huagui Duan$^{1}$,\thanks{Partially supported by NSFC (No.11131004, 11471169), LPMC of MOE of China and
Nankai University. E-mail: duanhg@nankai.edu.cn} \qquad
Hui Liu$^{2}$\thanks{Partially supported by NSFC (No.11401555), China Postdoctoral Science Foundation No.2014T70589, CUSF(No. WK0010000037).
E-mail:huiliu@ustc.edu.cn } \\ \\
$^{1}$ School of Mathematical Sciences and LPMC, Nankai University,\\ Tianjin 300071, P. R. China\\
$^{2}$ Wu Wen-Tsun Key Laboratory of Mathematics, USTC, Chinese Academy of Sciences,\\
School of Mathematical Sciences, University of Science and Technology of China,
\\Hefei, Anhui 230026, P. R. China\\}

\begin{document}
\date{April 21, 2015}
\maketitle

\begin{abstract}
{\it If all prime closed geodesics on $(S^n,F)$ with an irreversible Finsler metric $F$ are irrationally elliptic,
there exist either exactly $2\left[\frac{n+1}{2}\right]$ or infinitely many distinct closed geodesics. As an application,
we show the existence of three distinct closed geodesics on bumpy Finsler $(S^3,F)$ if any prime closed geodesic has non-zero Morse index.}
\end{abstract}



\renewcommand{\theequation}{\thesection.\arabic{equation}}
\renewcommand{\thefigure}{\thesection.\arabic{figure}}

\setcounter{equation}{0}
\section{Introduction and main results}%{Section 1}

There is a famous conjecture in Riemannian geometry which claims that there exist infinitely many
closed geodesics on any compact Riemannian manifold. This conjecture
has been proved except for CROSS's (compact rank one symmetric
spaces). The results of Franks \cite{Fra} in 1992 and Bangert \cite{Ban}
in 1993 imply that this conjecture is true for any Riemannian 2-sphere.
But once one move to the Finsler case, the conjecture
becomes false. It was quite surprising when Katok \cite{Kat}  in 1973 found
some irreversible Finsler metrics on CROSS's with only finitely
many closed geodesics and all closed geodesics being non-degenerate and elliptic.

{\bf Definition 1.1.} (cf. \cite{She}) {\it Let $M$ be a finite
dimensional smooth manifold. A function $F:TM\to [0,+\infty)$ is a {\rm
Finsler metric} if it satisfies

(F1) $F$ is $C^{\infty}$ on $TM\bs\{0\}$,

(F2) $F(x,\lm y) = \lm F(x,y)$ for all $x\in M$, $y\in T_xM$ and
$\lm>0$,

(F3) For every $y\in T_xM\bs\{0\}$, the quadratic form
$$ g_{x,y}(u,v) \equiv
         \frac{1}{2}\frac{\pt^2}{\pt s\pt t}F^2(x,y+su+tv)|_{t=s=0},
         \qquad \forall u, v\in T_xM, $$
is positive definite.

In this case, $(M,F)$ is called a {\rm Finsler manifold}. $F$ is
{\rm reversible} if $F(x,-y)=F(x,y)$ holds for all $x\in M$ and
$y\in T_xM$. $F$ is {\rm Riemannian} if $F(x,y)^2=\frac{1}{2}G(x)y\cdot
y$ for some symmetric positive definite matrix function $G(x)\in
GL(T_xM)$ depending on $x\in M$ smoothly. }

A closed curve on a Finsler manifold is a closed geodesic if it is
locally the shortest path connecting any two nearby points on this
curve (cf. \cite{She}). As usual, on any Finsler manifold
$(M, F)$, a closed geodesic $c:S^1=\R/\Z\to M$ is {\it prime}
if it is not a multiple covering (i.e., iteration) of any other
closed geodesics. Here the $m$-th iteration $c^m$ of $c$ is defined
by $c^m(t)=c(mt)$. The inverse curve $c^{-1}$ of $c$ is defined by
$c^{-1}(t)=c(1-t)$ for $t\in \R$.  Note that unlike Riemannian manifold,
the inverse curve $c^{-1}$ of a closed geodesic $c$
on a irreversible Finsler manifold need not be a geodesic.
We call two prime closed geodesics
$c$ and $d$ {\it distinct} if there is no $\th\in (0,1)$ such that
$c(t)=d(t+\th)$ for all $t\in\R$.
On a reversible Finsler (or Riemannian) manifold, two closed geodesics
$c$ and $d$ are called { \it geometrically distinct} if $
c(S^1)\neq d(S^1)$, i.e., their image sets in $M$ are distinct.
We shall omit the word {\it distinct} when we talk about more than one prime closed geodesic.

For a closed geodesic $c$ on $(S^n,\,F)$, denote by $P_c$
the linearized Poincar\'{e} map of $c$. Then $P_c\in \Sp(2n-2)$ is symplectic.
For any $M\in \Sp(2k)$, we define the {\it elliptic height } $e(M)$
of $M$ to be the total algebraic multiplicity of all eigenvalues of
$M$ on the unit circle $\U=\{z\in\C|\; |z|=1\}$ in the complex plane
$\C$. Since $M$ is symplectic, $e(M)$ is even and $0\le e(M)\le 2k$.
A closed geodesic $c$ is called {\it elliptic} if $e(P_c)=2(n-1)$, i.e., all the
eigenvalues of $P_c$ locate on $\U$; {\it irrationally elliptic} if, in the homotopy component $\Om^0(P_c)$
of $P_c$ (cf. Section 2 below), $P_c$ can be connected to the $\dm$-product of $(n-1)$ rotation matrices
$\left(\matrix{\cos\th_i & -\sin\th_i \cr \sin\th_i & \cos\th_i\cr}\right)$
with $\th_i$ being irrational multiples of $\pi$ for $1\le i\le n-1$; {\it hyperbolic} if $e(P_c)=0$, i.e., all the
eigenvalues of $P_c$ locate away from $\U$;
{\it non-degenerate} if $1$ is not an eigenvalue of $P_c$. A Finsler manifold $(M,\,F)$
is called {\it bumpy} if all the closed geodesics on it are non-degenerate.

Recently the closed geodesic problem on Finsler spheres has aroused many authors' interest
(cf. \cite{BaL}, \cite{DuL1}-\cite{DuL3}, \cite{HWZ1}-\cite{HWZ2}, \cite{Lon5}, \cite{LoD}, \cite{LoW}, \cite{Rad3}-\cite{Rad4}, \cite{Wan} and so on).

Note that Hingston in \cite{Hin} proved the existence of infinitely many closed geodesics on Riemannian spheres
if all prime closed geodesics are hyperbolic. The results of \cite{HWZ1}-\cite{HWZ2} imply that there exist either
two or infinitely many prime closed geodesics on every bumpy Finsler 2-sphere provided the stable and unstable manifolds
of every hyperbolic closed geodesic intersect transversally. After the work \cite{BaL} of Bangert and Long, Long and Wang
in \cite{LoW} proved that on every Finsler 2-sphere, there exist either at least two irrationally elliptic closed geodesics
or infinitely many closed geodesics.

In the Katok's examples of \cite{Kat} which have been studied later by Ziller in \cite{Zil}, Katok constructed
some irreversible Finsler metrics  on $S^n$ with exact $2\left[\frac{n+1}{2}\right]$ closed geodesics and all
of them being irrationally elliptic.  Based on these facts, Anosov in \cite{Ano} conjectured that there exist
at least $2\left[\frac{n+1}{2}\right]$ closed geodesics on every Finsler spheres $(S^n,F)$.
Furthermore one suspects that all these closed geodesics on $S^n$ with an irreversible Finsler metric $F$ are irrationally
elliptic provided the number of prime closed geodesics is finite (cf. \cite{Lon5} for more conjectures of the closed geodesic problem).
See \cite{WHL} for a similar conjecture about closed characteristics on compact convex hypersurfaces in $\R^{2n}$.

Thus it seems interesting to study the multiplicity problem on the Finsler sphere $(S^n,F)$ with irrationally elliptic closed geodesics.
The following is our main result.

\medskip

{\bf Theorem 1.2.} {\it If all prime closed geodesics on $(S^n,F)$ with an irreversible Finsler metric $F$ are irrationally elliptic,
there exist either exactly $2\left[\frac{n+1}{2}\right]$ or infinitely many distinct closed geodesics.}

\medskip

In \cite{DuL2}, Duan and Long showed that on every bumpy Finsler $(S^3,F)$ either there exist at least three distinct closed geodesics,
or there exists two non-hyperbolic closed geodesics. As an application of Theorem 1.2, we obtain

\medskip

{\bf Theorem 1.3.}  {\it There exist at least three distinct closed geodesics on a bumpy Finsler $(S^3,F)$ if any prime closed geodesic on it has non-zero Morse index.}

\medskip

In this paper, let $\N$, $\N_0$, $\Z$, $\Q$, $\R$, and $\C$ denote
the sets of natural integers, non-negative integers, integers,
rational numbers, real numbers, and complex numbers respectively.
We use only singular homology modules with $\Q$-coefficients.
For an $S^1$-space $X$, we denote by $\overline{X}$ the quotient space $X/S^1$.
We define the functions
\be \left\{\matrix{[a]=\max\{k\in\Z\,|\,k\le a\}, &
           E(a)=\min\{k\in\Z\,|\,k\ge a\} , \cr
    \varphi(a)=E(a)-[a], &\{a\}=a-[a].  \cr}\right. \lb{1.1}\ee
Especially, $\varphi(a)=0$ if $ a\in\Z\,$, and $\varphi(a)=1$ if $
a\notin\Z\,$.

\setcounter{equation}{0}
\section{Morse theory and Morse indices of closed geodesics}%{Section 2}

Let $M=(M,F)$ be a compact Finsler manifold $(M,F)$, the space
$\Lambda=\Lambda M$ of $H^1$-maps $\gamma:S^1\rightarrow M$ has a
natural structure of Riemannian Hilbert manifolds on which the
group $S^1=\R/\Z$ acts continuously by isometries (cf. \cite{Kli}). This action is defined by
$(s\cdot\gamma)(t)=\gamma(t+s)$ for all $\gamma\in\Lm$ and $s,
t\in S^1$. For any $\gamma\in\Lambda$, the energy functional is
defined by
\be E(\gamma)=\frac{1}{2}\int_{S^1}F(\gamma(t),\dot{\gamma}(t))^2dt.
\lb{2.1}\ee
It is $C^{1,1}$ and invariant under the $S^1$-action. The
critical points of $E$ of positive energies are precisely the closed geodesics
$\gamma:S^1\to M$. The index form of the functional $E$ is well
defined along any closed geodesic $c$ on $M$, which we denote by
$E''(c)$. As usual, we denote by $i(c)$ and
$\nu(c)$ the Morse index and nullity of $E$ at $c$. In the
following, we denote by
\be \Lm^\kappa=\{d\in \Lm\;|\;E(d)\le\kappa\},\quad \Lm^{\kappa-}=\{d\in \Lm\;|\; E(d)<\kappa\},
  \quad \forall \kappa\ge 0. \nn\ee
For a closed geodesic $c$ we set $ \Lm(c)=\{\ga\in\Lm\mid E(\ga)<E(c)\}$.

For $m\in\N$ we denote the $m$-fold iteration map
$\phi_m:\Lambda\rightarrow\Lambda$ by $\phi_m(\ga)(t)=\ga(mt)$, for all
$\,\ga\in\Lm, t\in S^1$, as well as $\ga^m=\phi_m(\gamma)$. If $\gamma\in\Lambda$
is not constant then the multiplicity $m(\gamma)$ of $\gamma$ is the order of the
isotropy group $\{s\in S^1\mid s\cdot\gamma=\gamma\}$. For a closed geodesic $c$,
the mean index $\hat{i}(c)$ is defined as usual by
$\hat{i}(c)=\lim_{m\to\infty}i(c^m)/m$. Using singular homology with rational
coefficients we consider the following critical $\Q$-module of a closed geodesic
$c\in\Lambda$:
\be \overline{C}_*(E,c)
   = H_*\left((\Lm(c)\cup S^1\cdot c)/S^1,\Lm(c)/S^1\right). \lb{2.3}\ee

{\bf Proposition 2.1.} (cf. Satz 6.11 of \cite{Rad2} ) {\it Let $c$ be a
prime closed geodesic on a bumpy Finsler manifold $(M,F)$. Then there holds}
$$ \overline{C}_q( E,c^m) = \left\{\matrix{
     \Q, &\quad {\it if}\;\; i(c^m)-i(c)\in 2\Z\;\;{\it and}\;\;
                   q=i(c^m),\;  \cr
     0, &\quad {\it otherwise}. \cr}\right.  $$

Set $\ol{\Lm}^0=\ol{\Lambda}^0M =\{{\rm
constant\;point\;curves\;in\;}M\}\cong M$. Let $(X,Y)$ be a
space pair such that the Betti numbers $b_i=b_i(X,Y)=\dim
H_i(X,Y;\Q)$ are finite for all $i\in \Z$. As usual the {\it
Poincar\'e series} of $(X,Y)$ is defined by the formal power series
$P(X, Y)=\sum_{i=0}^{\infty}b_it^i$. We need the following well
known version of results on Betti numbers.

{\bf Lemma 2.2.} (cf. Theorem 2.4 and Remark 2.5 of \cite{Rad1} and \cite{Hin}, cf.
also Lemma 2.5 of \cite{DuL3}) {\it Let $(S^n,F)$ be an
$n$-dimensional Finsler sphere.}

(i) {\it When $n$ is odd, the Betti numbers are given by
\bea b_j
&=& \rank H_j(\Lm S^n/S^1,\Lm^0 S^n/S^1;\Q)  \nn\\
&=& \left\{\matrix{
    2,&\quad {\it if}\quad j\in \K\equiv \{k(n-1)\,|\,2\le k\in\N\},  \cr
    1,&\quad {\it if}\quad j\in \{n-1+2k\,|\,k\in\N_0\}\bs\K,  \cr
    0 &\quad {\it otherwise}. \cr}\right. \lb{2.4}\eea}

(ii) {\it When $n$ is even, the Betti numbers are given by
\bea b_j
&=& \rank H_j(\Lm S^n/S^1,\Lm^0 S^n/S^1;\Q)  \nn\\
&=& \left\{\matrix{
    2,&\quad {\it if}\quad j\in \K\equiv \{k(n-1)\,|\,3\le k\in 2\N+1\},  \cr
    1,&\quad {\it if}\quad j\in \{n-1+2k\,|\,k\in\N_0\}\bs\K,  \cr
    0 &\quad {\it otherwise}. \cr}\right.    \lb{2.5}\eea}

Next we recall some index iterated formulae of closed geodesics (cf. \cite{Bot} and \cite{Lon1}-\cite{Lon4}).

In \cite{Lon2} of 1999, Y. Long established the basic normal form
decomposition of symplectic matrices. Based on this result he
further established the precise iteration formulae of indices of
symplectic paths in \cite{Lon3} of 2000. Note that this index iteration formulae works for Morse indices
of iterated closed geodesics (cf. \cite{LiL}, \cite{Liu} and Chap. 12 of \cite{Lon4}).
Since every closed geodesic on a sphere is orientable. Then by Theorem 1.1 of \cite{Liu}, the initial Morse index of a closed geodesic
$c$ on a $n$-dimensional Finsler  sphere coincides with the index of a
corresponding symplectic path.

As in \cite{Lon4}, denote by
\bea
N_1(\lm, a) &=& \left(\matrix{\lm & a\cr
                                0 & \lm\cr}\right), \qquad {\rm for\;}\lm=\pm 1, \; a\in\R, \lb{2.6}\\
H(b) &=& \left(\matrix{b & 0\cr
                      0 & b^{-1}\cr}\right), \qquad {\rm for\;}b\in\R\bs\{0, \pm 1\}, \lb{2.7}\\
R(\th) &=& \left(\matrix{\cos\th & -\sin\th \cr
                           \sin\th & \cos\th\cr}\right), \qquad {\rm for\;}\th\in (0,\pi)\cup (\pi,2\pi), \lb{2.8}\\
N_2(e^{\th\sqrt{-1}}, B) &=& \left(\matrix{ R(\th) & B \cr
                  0 & R(\th)\cr}\right), \qquad {\rm for\;}\th\in (0,\pi)\cup (\pi,2\pi)\;\; {\rm and}\; \nn\\
        && \qquad B=\left(\matrix{b_1 & b_2\cr
                                  b_3 & b_4\cr}\right)\; {\rm with}\; b_j\in\R, \;\;
                                         {\rm and}\;\; b_2\not= b_3. \lb{2.9}\eea
Here $N_2(e^{\th\sqrt{-1}}, B)$ is non-trivial if $(b_2-b_3)\sin\theta<0$, and trivial
if $(b_2-b_3)\sin\theta>0$ as defined in \cite{Lon3} and Definition 1.8.11 of \cite{Lon4}.

As in \cite{Lon4}, the $\diamond$-sum (direct sum) of any two real matrices is defined by
$$ \left(\matrix{A_1 & B_1\cr C_1 & D_1\cr}\right)_{2i\times 2i}\diamond \left(\matrix{A_2 & B_2\cr C_2 & D_2\cr}\right)_{2j\times 2j}
=\left(\matrix{A_1 & 0 & B_1 & 0 \cr
                                   0 & A_2 & 0& B_2\cr
                                   C_1 & 0 & D_1 & 0 \cr
                                   0 & C_2 & 0 & D_2}\right). $$

For every
$P\in\Sp(2n)$, the homotopy set $\Omega(P)$ of $P$ in $\Sp(2n)$ is defined by
$$ \Om(P)=\{N\in\Sp(2n)\,|\,\sg(N)\cap\U=\sg(P)\cap\U\equiv\Gamma\;\mbox{and}
                    \;\nu_{\om}(N)=\nu_{\om}(P)\, \forall\om\in\Gamma\}, $$
where $\sg(P)$ denotes the spectrum of $P$,
$\nu_{\om}(P)\equiv\dim_{\C}\ker_{\C}(P-\om I)$ for $\om\in\U$.
The homotopy component $\Om^0(P)$ of $P$ in $\Sp(2n)$ is defined by
the path connected component of $\Om(P)$ containing $P$. Then the following decomposition theorem is proved in \cite{Lon2}
and \cite{Lon3}.

\medskip

{\bf Theorem 2.3.} (cf. Theorem 1.8.10, Lemma 2.3.5 and Theorem 8.3.1 of \cite{Lon4}) {\it  For every $P\in\Sp(2n)$, there
exists a continuous path $f\in\Om^0(P)$ such that $f(0)=P$ and
\bea f(1)
&=& N_1(1,1)^{\dm p_-}\,\dm\,I_{2p_0}\,\dm\,N_1(1,-1)^{\dm p_+}
  \dm\,N_1(-1,1)^{\dm q_-}\,\dm\,(-I_{2q_0})\,\dm\,N_1(-1,-1)^{\dm q_+} \nn\\
&&\dm\,N_2(e^{\aa_{1}\sqrt{-1}},A_{1})\,\dm\,\cdots\,\dm\,N_2(e^{\aa_{r_{\ast}}\sqrt{-1}},A_{r_{\ast}})
  \dm\,N_2(e^{\bb_{1}\sqrt{-1}},B_{1})\,\dm\,\cdots\,\dm\,N_2(e^{\bb_{r_{0}}\sqrt{-1}},B_{r_{0}})\nn\\
&&\dm\,R(\th_1)\,\dm\,\cdots\,\dm\,R(\th_k)\,\dm\,R(\th_{k+1})\,\dm\,\cdots\,\dm\,R(\th_r)\dm\,H(b)^{\dm h},\lb{2.10}\eea
where $\frac{\th_{j}}{2\pi}\not\in\Q\cap(0,1)$ for $1\le j\le k$ and
$\frac{\th_{j}}{2\pi}\in\Q\cap(0,1)$ for $k+1\le j\le r$; $N_2(e^{\aa_{j}\sqrt{-1}},A_{j})$'s
are nontrivial and $N_2(e^{\bb_{j}\sqrt{-1}},B_{j})$'s are trivial, and non-negative integers
$p_-, p_0, p_+,q_-, q_0, q_+,r,r_\ast,r_0,h$ satisfy
\be p_- + p_0 + p_+ + q_- + q_0 + q_+ + r + 2r_{\ast} + 2r_0 + h = n. \lb{2.11}\ee

Let $\ga\in\P_{\tau}(2n)=\{\ga\in C([0,\tau],\Sp(2n))\,|\,\ga(0)=I\}$. Denote the basic normal form
decomposition of $P\equiv \ga(\tau)$ by (\ref{2.10}). Then we have
\bea i(\ga^m)
&=& m(i(\ga)+p_-+p_0-r ) + 2\sum_{j=1}^rE\left(\frac{m\th_j}{2\pi}\right) - r   \nn\\
&&  - p_- - p_0 - {{1+(-1)^m}\over 2}(q_0+q_+)
              + 2\sum_{j=1}^{r_{\ast}}\vf\left(\frac{m\aa_j}{2\pi}\right) - 2r_{\ast}. \lb{2.12}\eea}

\setcounter{equation}{0}
\section{Proof of main theorems}

At first we prove Theorem 1.2. Now we assume the following condition holds,

{\bf (IECG)}: Suppose that all prime closed geodesics on $(S^n,F)$ with finitely many distinct closed geodesics are irrationally elliptic, denoted by $\{c_j\}_{j=1}^q$.

{\bf Proof of Theorem 1.2:}

Under the assumption ({\bf IECG}), by the decomposition in Theorem 2.3, the linearized Poincar\'{e} map $P_{c_j}$
can be connected to $f_{c_j}(1)$ in $\Om^0(P_{c_j})$ satisfying
\bea f_{c_j}(1)&=&R(\th_{j1})\,\dm\,\cdots\,\dm\,R(\th_{j(n-1)}), \quad 1\le j\le n-1, \lb{3.1}\eea
 where $\frac{\th_{jk}}{2\pi}\not\in\Q$ for $1\le j\le q$ and $1\le k\le n-1$.  Then by Theorem 2.3, we obtain
\bea i(c_j)&\in& 2\Z+(n-1),\qquad 1\le j\le q,\lb{3.2}\\
i(c_j^m)&=& m(i(c_j)-n+1) + 2\sum_{k=1}^{n-1}\left[\frac{m\th_{jk}}{2\pi}\right] +n-1,\quad \nu(c_j^m)=0,\qquad\forall\ m\ge 1, \lb{3.3}\eea
where (\ref{3.2}) follows from Theorem 8.1.7 of \cite{Lon4} and the symplectic additivity of iterated indices.

By (\ref{3.2}) and (\ref{3.3}), there holds
\bea i(c_j^{m+1})-i(c_j^m)\in 2\Z,\qquad \forall\ 1\le j\le q,\quad m\ge 1. \lb{3.4}\eea

Define
\bea M_p=\sum_{j=1}^q M_p(j)\equiv\sum_{j=1}^q\#\{m\ge 1\ |\ i(c_j^m)=p, \;\ol{C}_p(E, c_j^m)\not= 0\},\quad p\in\Z.\nn\eea
Then the following Morse inequality (cf. Theorem I.4.3 of \cite{Cha}) holds
\bea M_p&\ge& b_p,\lb{3.5}\\
M_p - M_{p-1} + \cdots +(-1)^{p}M_0
&\ge& b_p - b_{p-1}+ \cdots + (-1)^{p}b_0, \quad\forall\ p\in\N_0.\lb{3.6}\eea

By Proposition 2.1 and (\ref{3.4}) we obtain
\bea M_p=\sum_{j=1}^q M_p(j)=\sum_{j=1}^q\#\{m\ge 1\ |\ i(c_j^m)=p\},\qquad \forall\ p\in\Z, \lb{3.7}\eea
which, together with (\ref{3.2}) and Lemma 2.2, yields
\be M_{p}=b_{p}=0, \qquad \forall\ p=n\ (\mod 2),\ p\in\N_0.{\lb{3.8}}\ee

Then by the Morse inequality (\ref{3.6}), we obtain
\be M_{p}=b_{p}, \qquad \forall\ p=n-1\ (\mod 2),\ p\in\N_0.{\lb{3.9}}\ee

{\bf Claim 1:} {\it Under the condition (IECG), there holds $i(c_j)\ge n-1$, $\forall\ 1\le j\le q$. Moreover, $\#\{1\le j\le q\ |\ i(c_j)=n-1\}=1$.}

\medskip

In fact, if there exists at least one closed geodesic $c_{j_0}$ such that $0\leq i(c_{j_0})< n-1$, then by the definition of $M_p$ it yields
$M_{i(c_{j_0})}\ge 1$. Note that $i(c_{j_0})=n-1\ (\mod 2)$ by (\ref{3.2}). Then by the Morse inequality and Lemma 2.2 we can obtain the following contradiction
\be -1\ge M_{i(c_{j_0})+1} - M_{i(c_{j_0})} + \cdots +(-1)^{i(c_{j_0})+1}M_0
\ge b_{i(c_{j_0})+1} - b_{i(c_{j_0})}+ \cdots + (-1)^{i(c_{j_0})+1}b_0=0.\lb{3.10}\ee

In addition, by (\ref{3.9}) and Lemma 2.2, there holds $M_{n-1}=b_{n-1}=1$. So $\#\{1\le j\le q\ |\ i(c_j)=n-1\}=1$ by (\ref{3.7}). So Claim 1 holds.

\medskip

By Claim 1, without loss of generality, assume that $i(c_1)=n-1$ and $i(c_j)\ge n+1$, $\forall\ 2\le j\le q$. Then by (\ref{3.3}), it yields $\hat{i}(c_j)=i(c_j)-(n-1)+
\sum_{k=1}^{n-1}\frac{\th_{jk}}{\pi}>0$, $\forall\ 1\le j\le q$. Thus by the common index jump theorem (cf. Theorems 4.1-4.3 of \cite{LoZ}), there exist infinitely many
$(N, m_1,\ldots,m_q)\in\N^{q+1}$ such that
\bea
i(c_1^{2m_1-1}) &=& 2N-i(c_1)=2N-(n-1), \lb{3.11}\\
i(c_1^{2m_1+1}) &=& 2N+i(c_1)=2N+(n-1), \lb{3.12}\\
2N-(n-1)&\le& i(c_1^{2m_1}) \le 2N+(n-1), \lb{3.13}\\
i(c_j^{2m_j-1}) &=& 2N-i(c_j)\le 2N-(n+1), \lb{3.14}\\
i(c_j^{2m_j+1}) &=& 2N+i(c_j)\ge 2N+(n+1), \lb{3.15}\\
2N-(n-1)&\le& i(c_j^{2m_j}) \le 2N+(n-1), \quad 2\le j\le q,\lb{3.16}
\eea
where, furthermore, given $M_0\in\N$, we can require  $M_0|N$ (Note that the closure of the
set $\{\{Nv\}: N\in\N, \;M_0|N\}$ is still a closed additive subgroup of $\bf T^h$ for some $h\in\N$,
where $v$ is defined in (4.21) of \cite{LoZ}. Then we can use the proof of Step 2 in Theorem 4.1
of \cite{LoZ} to get $N$). The freedom of choice of $N$ will be used in the proof of Claim 3 below.

\medskip

{\bf Claim 2:}  {\it There holds $i(c_1^{2m_1+m})\ge 2N+(n+1)$, $i(c_1^{2m_1-m})\le 2N-(n+1)$ for any $m\ge2$.}

\medskip

In fact, by $\#\{m\ge 1\ |\ i(c_1^m)=n-1\}=M_{n-1}=b_{n-1}=1$ we get $i(c_1^m)\ge n+1,\ \forall\ m\ge 2$. So by (\ref{3.3}) and $i(c)=n-1$ we have
\bea i(c_1^2)&=& 2\sum_{k=1}^{n-1}\left[\frac{\th_{1k}}{\pi}\right] +n-1\ge n+1,\eea which implies that there exists at least some rotation angel $\th_{1r}$ with $1\le r\le n-1$ such that $1<\frac{\th_{1r}}{\pi}<2$. Therefore again by (\ref{3.3}) and $i(c_1)=n-1$ we obtain
\bea i(c_1^{2m_1+2})-i(c_1^{2m_1+1})&=& 2\sum_{k=1}^{n-1}\left(\left[\frac{m_1\th_{1k}}{\pi}+\frac{\th_{1k}}{\pi}\right]-\left[\frac{m_1\th_{1k}}{\pi}+\frac{\th_{1k}}{2\pi}\right]\right)\nn\\
         &\ge&2\left(\left[\frac{m_1\th_{1r}}{\pi}+\frac{\th_{1r}}{\pi}\right]-\left[\frac{m_1\th_{1r}}{\pi}+\frac{\th_{1r}}{2\pi}\right]\right)\nn\\
         &=&2\left(\left[\left\{\frac{m_1\th_{1r}}{\pi}\right\}+\left\{\frac{\th_{1r}}{\pi}\right\}\right]-\left[\left\{\frac{m_1\th_{1r}}{\pi}\right\}+\frac{\th_{1r}}{2\pi}\right]\right)+2\nn\\
         &=&2,\lb{3.18}\eea
where the last equality holds by choosing
\bea \left\{\frac{m_1\th_{1r}}{\pi}\right\}>\max\left\{1-\left\{\frac{\th_{1r}}{\pi}\right\}, 1-\frac{\th_{1r}}{2\pi}\right\}.\lb{3.19}\eea

Similarly we can obtain
\bea i(c_1^{2m_1-1})-i(c_1^{2m_1-2})&=& 2\sum_{k=1}^{n-1}\left(\left[\frac{m_1\th_{1k}}{\pi}-\frac{\th_{1k}}{2\pi}\right]-\left[\frac{m_1\th_{1k}}{\pi}-\frac{\th_{1k}}{\pi}\right]\right)\nn\\
         &\ge&2\left(\left[\frac{m_1\th_{1r}}{\pi}-\frac{\th_{1r}}{2\pi}\right]-\left[\frac{m_1\th_{1r}}{\pi}-\frac{\th_{1r}}{\pi}\right]\right)\nn\\
         &=&2\left(\left[\left\{\frac{m_1\th_{1r}}{\pi}\right\}-\frac{\th_{1r}}{2\pi}\right]-\left[\left\{\frac{m_1\th_{1r}}{\pi}\right\}-\left\{\frac{\th_{1r}}{\pi}\right\}\right]\right)+2\nn\\
         &=&2,\lb{3.20}\eea
where the last equality holds by choosing
\bea \left\{\frac{m_1\th_{1r}}{\pi}\right\}>\max\left\{\left\{\frac{\th_{1r}}{\pi}\right\}, \frac{\th_{1r}}{2\pi}\right\}.\lb{3.21}\eea

Here the existence of $m_1$ satisfying requirements (\ref{3.19}) and (\ref{3.21}) can be made sure by (4.5) of Theorem 4.1, (4.43) in the proof of Theorem 4.3 and (c) of Theorem 4.2 in \cite{LoZ}, or directly by Corollary 3.19 of \cite{DuL3}.

In addition, it follows from $i(c_1)=n-1$ and (\ref{3.3}) that $i(c_1^{m+1})\ge i(c_1^m),\ \forall\ m\ge1$. Then by (\ref{3.11})-(\ref{3.12}), (\ref{3.18}) and (\ref{3.20}), Claim 2 is proved.

\medskip

{\bf Claim 3:} {$q = 2\left[\frac{n+1}{2}\right]$.}

\medskip

We use two cases to carry out the proof of Claim 3:

{\bf Case 1:} {\it $n\in 2\N$.}

For the integer $N$ in the common index jump theorem, we can require $(n-1)|N$. Thus $b_{2N+n-1}=b_{2N-n+1}=2, b_p=1, \forall\ p\in (2N-(n-1), 2N+(n-1))\cap (2\N-1)$
by Lemma 2.2. Thus together with (\ref{2.5}),
it yields \bea \sum_{p=2N-(n-1)}^{2N+n-1}b_p=\frac{(2N+n-1)-(2N-n+1)}{2}+3=n+2.\lb{3.22}\eea

On the other hand, by (\ref{3.7}), (\ref{3.11})-(\ref{3.16}) and Claim 2 we have
\bea \sum_{p=2N-(n-1)}^{2N+n-1}M_p(1)=3, \qquad \sum_{p=2N-(n-1)}^{2N+n-1}M_p(j)=1,\quad \forall\ 2\le j\le q,\nn\eea
which yields
\bea \sum_{p=2N-(n-1)}^{2N+n-1}M_p=3+q-1=q+2.\lb{3.23}\eea

Thus from (\ref{3.8})-(\ref{3.9}), (\ref{3.22}) and (\ref{3.23}) it follows
\bea n+2=\sum_{p=2N-(n-1)}^{2N+n-1}b_p=\sum_{p=2N-(n-1)}^{2N+n-1}M_p=q+2,\eea
which yields $q=n$ and completes the proof of Claim 3 in Case 1.

{\bf Case 2:} {\it $n\in 2\N-1$.}

For the integer $N$ in the common index jump theorem, we can require $(n-1)|N$. Thus $b_{2N+n-1}=b_{2N}=b_{2N-n+1}=2, b_p=1, \forall\ p\in 2\N\cap(2N-(n-1), 2N+(n-1))\backslash\{2N\}$ by Lemma 2.2. Thus similarly to the proof in Case 1
it yields \bea n+3&=&\frac{(2N+n-1)-(2N-n+1)}{2}+4\nn\\
     &=&\sum_{p=2N-(n-1)}^{2N+n-1}b_p\nn\\
     &=&\sum_{p=2N-(n-1)}^{2N+n-1}M_p\nn\\
     &=&q+2,\lb{3.25}\eea
which yields $q=n+1$ and completes the proof of Claim 3 in Case 2.

By the assumption ({\bf IECG}) and Claim 3, we finish the proof of Theorem 1.2.\hfill\hb

\medskip

{\bf Proof of Theorem 1.3:}

\medskip

It is proved in Remark 3.7 of \cite{DuL2} that suppose that there exist precisely two prime closed geodesics $c_1$ and $c_2$ on every bumpy Finsler $(S^3,F)$,
then both $c_1$ and $c_2$ are non-hyperbolic and must belong to one of three precise classes there. Furthermore, if any prime closed geodesic on it has
non-zero Morse index, then both $c_1$ and $c_2$ are irrationally elliptic. Thus as an application of Theorem 1.2 we proved Theorem 1.3. \hfill\hb

\medskip

{\bf Acknowledgements.} The authors would like to thank sincerely Professor Yiming Long for his valuable help and encouragement to them. And the authors also thank sincerely him for his comments, suggestions and helpful discussions about the closed geodesic problem with them.

\bibliographystyle{abbrv}

\end{document}